\documentclass[reqno]{amsart}
\usepackage{amssymb}
\usepackage{amsmath}

\makeatletter
\@addtoreset{equation}{section}
\makeatother

\renewcommand\thefigure{\thesection.\@arabic\c@figure}
\renewcommand\thetable{\thesection.\@arabic\c@table}

\newtheorem{theorem}{Theorem}[section]
\newtheorem{lemma}[theorem]{Lemma}
\newtheorem{proposition}[theorem]{Proposition}
\newtheorem{corollary}[theorem]{Corollary}

\newcommand{\mc}[1]{{\mathcal #1}}
\newcommand{\mf}[1]{{\mathfrak #1}}
\newcommand{\mb}[1]{{\mathbf #1}}
\newcommand{\bb}[1]{{\mathbb #1}}

\newcommand{\<}{\langle}
\renewcommand{\>}{\rangle}

\newcommand{\x}{{\bf x}}
\def\Z{\mathbb Z}
\def\R{\mathbb R}
\def\N{\mathbb N}
\def\V{\mathcal V}

\begin{document}

\author{M. D. Jara, C. Landim, and S. Sethuraman}

\address{\noindent IMPA, Estrada Dona Castorina 110,
CEP 22460 Rio de Janeiro, Brasil
\newline
e-mail:  \rm \texttt{monets@impa.br}
}

\address{\noindent IMPA, Estrada Dona Castorina 110,
CEP 22460, Rio de Janeiro, Brasil and CNRS UMR 6085,
Avenue de l'Universit\'e, BP.12, Technop\^ole du Madrillet,
F76801 Saint-\'Etienne-du-Rouvray, France.
\newline
e-mail:  \rm \texttt{landim@impa.br}
}

\address{\noindent Department of Mathematics,
396 Carver Hall, Iowa State University, Ames, IA \ 50011, USA.
\newline
e-mail:  \rm \texttt{sethuram@iastate.edu}}

\title[Nonequilibrium fluctuations for a tagged particle]
{Nonequilibrium fluctuations for a tagged particle in mean-zero
one-dimensional zero-range processes}

\begin{abstract}
  We prove a non-equilibrium functional central limit theorem for the position of
  a tagged particle in mean-zero one-dimensional zero-range
  process. The asymptotic behavior of the tagged particle is described
  by a stochastic differential equation governed by the solution of
  the hydrodynamic equation.
\end{abstract}

\subjclass[2000]{primary 60K35}

\keywords{Hydrodynamic limit, tagged particle,
scaling limit, nonequilibrium}

\thanks{C. Landim was partially supported by the John S. Guggenheim
  Memorial Foundation, FAPERJ and CNPq; M.D. Jara and S. Sethuraman were
  partially supported by NSA-H982300510041 and NSF-DMS-0504193.}

\maketitle

\section{Introduction}
\label{sec0}

Informally, the zero-range particle system follows a collection of
dependent random walks on the lattice $\Z^d$ where, from a vertex with
$k$ particles, one of the particles displaces by $j$ with rate
$(g(k)/k)p(j)$.  The function on the non-negative integers $g:\N_0
\rightarrow \R_+$ is called the process ``rate'', and $p(\cdot)$
denotes the translation-invariant single particle transition
probability. The name ``zero-range'' derives from the observation
that, infinitesimally, the interaction is only with respect those
particles at the particular vertex.  The case when $g(k)$ is
proportional to $k$ describes the situation of completely independent
particles.
\medskip

The problem of the asymptotics of a distinguished, or tagged particle
interacting with others has a long history and was even mentioned in
Spitzer's paper \cite{Spitzer} (see also chapters 8.I, 6.II
\cite{Spohn}). The main analytical difficulty is that the tagged
particle motion is not in general Markovian due to the interaction
with other particles. However, the intuition is that in a scale the
tagged particle behaves as a random walk with certain ``homogenized''
parameters reflecting the system dynamics.


We prove in this article a nonequilibrium invariance principle, with
respect to a diffusion process whose coefficients depend on the
hydrodynamic density, for the diffusively rescaled position of the tagged
particle in one-dimensional zero-range processes when the transition
probability $p$ is finite-range and mean-zero.  This invariance
principle is the first result which captures the nonequlibrium
fluctuations of a single, given particle in a general finite-range
interacting particle system. We remark, however, in \cite{JL}, a
nonequilibrium central limit{ theorem was proved for a tagged
particle in the nearest-neighbor symmetric one-dimensional simple
exclusion model by completely different methods which rely on the
special structure of the nearest-neighbor one-dimensional dynamics.
Also, we note, in \cite{Reza-pr}, a ``propagation of chaos'' type
nonequilibrium result was shown for finite-range symmetric $d\geq 1$
dimensional simple exclusion processes which gives the fluctuations
for a tagged particle selected at random, or in other words the
average tagged particle position; however, this result, which makes
key use of the ``averaging,'' does not convey the fluctuations of
any fixed, given particle and so is weaker than the one we state in
this paper.

We mention also, with respect to zero-range tagged particles,
previous results on laws of large numbers, in equilibrium
\cite{Saada}, \cite{Sext} and non-equilibrium \cite{Reza-lln}, and
equilibrium central limit theorems when the jump probability $p$ is
mean-zero, $\sum j p(j)=0$ \cite{Saada}, \cite{Sext}, and also when
$p$ is totally asymmetric and nearest-neighbor in $d=1$
\cite{Szrtg}, and also some diffusive variance results when $p$ has
a drift $\sum jp(j)\neq 0$ in $d=1$ and $d\geq 3$ \cite{Szrtg}.

\medskip

Denote by $\xi \in \bb N_0^{\bb Z}$, $\bb N_0 = \{0, 1, \dots\}$, the
states of the zero-range process, so that $\xi(x)$, $x\in\bb Z$,
stands for the total number of particles at site $x$ for the
configuration $\xi$.

Fix an integer $N\ge 1$, scale space by $N^{-1}$ and assume that the
zero-range process rescaled diffusively, $\{\xi_t^N : t\ge 0\}$,
starts from a local equilibrium state with density profile $\rho_0 :
\bb R\to \bb R_+$. Denote by $\{\pi^{N,0}_t : t\ge 0\}$ its empirical
measure. It is well known that $\pi^{N,0}_t$ converges in probability
to the absolutely continuous measure $\rho(t,u) du$, where $\rho(t,u)$
is the solution of a non-linear parabolic equation with initial
condition $\rho_0$.

Tag a particle initially at the origin and denote by $X^N_t$ its
position at time $t$. It is relatively simple to show that the
rescaled trajectory $\{X^N_t/N : 0\le t\le T\}$ is tight for the
uniform topology. In particular, to prove convergence, one needs
only to characterize the limit points.

In contrast with other models, in zero-range processes $X^N_t$ is a
square integrable martingale with a bounded quadratic variation
$\<X^N\>_t$ given by the time integral of a local function of the
process as seen from the tagged particle:
\begin{equation*}
\<X^N\>_t \;=\; \sigma^2 N^2 \int_0^t \frac {g(\eta^N_s(0))}
{\eta^N_s(0)} \, ds\;,
\end{equation*}
where $\sigma^2$ is the variance of the transition probability
$p(\cdot)$, $g(\cdot)$ is the jump rate mentioned before, and $\eta^N_s = \tau_{X^N_s}
\xi^N_s$ is the state of the process as seen from the tagged
particle. Here $\{\tau_x : x\in\bb Z\}$ stands for the group of
translations. In particular, if the rescaled position of the tagged
particle $x^N_t = X^N_t/N$ converges to some path $x_t$, this process
$x_t$ inherits the martingale property from $X^N_t$. If in addition
$x_t$ is continuous, to complete the
characterization, one needs to examine
the asymptotic behavior of its quadratic variation.

Denote by $\{\nu_\rho : \rho\ge 0\}$ the one-parameter family, indexed
by the density, of invariant states for the process as seen from the
tagged particle. Let $\pi^N_t$ be the empirical measure associated to
this process: $\pi^N_t = \tau_{X^N_t} \pi^{N,0}_t$ and suppose that
one can replace the local function $g(\eta^N_s(0))/\eta^N_s(0)$ by a
function of the empirical measure.  If we assume conservation of local
equilibrium for the process as seen from the tagged particle, this
function should be $h(\lambda (s,0))$, where $h(\rho)$ is the
expected value of $g(\eta(0))/\eta(0)$ under the invariant state
$\nu_\rho$ and $\lambda(s,0)$ is the density of particles around the
tagged particle, i.e., the density of particles around the origin for
the system as seen from the tagged particle.

As we are assuming that $X^N_t/N$ converges to $x_t$, since
$\pi^N_t = \tau_{X^N_t} \pi^{N,0}_t$ and $\pi^{N,0}_t$ converges
to $\rho(t,u) du$, we must have $\lambda(s,0) =
\rho(s,x_s)$. Therefore, if the quadratic variation of $X^N_t/N$
converges to the quadratic variation of $x_t$, $\<x\>_t = \sigma^2
\int_0^t h(\rho(s,x_s)) ds$. In particular, by the characterization of
continuous martingales, $x_t$ satisfies the stochastic
differential equation
\begin{equation*}
  dx_t \; =\; \sigma \sqrt{h(\rho(s,x_s))} \, dB_s\;,
\end{equation*}
where $\rho$ is the solution of the hydrodynamic equation, $h$ is
defined above and $B$ is a Brownian motion.

We see from this sketch that the main difficulty consists in proving
the conservation of local equilibrium around the tagged particle,
without assuming any type of attractiveness, which is relied upon in \cite{Landim-conservation}.  The absence of a space
average creates a major obstacle in this step. In contrast with the
proof of the hydrodynamic limit, we need to replace a local function
instead of a space average of translations of a local function.  We
may, therefore, only use the bonds close to the origin of the
Dirichlet form to perform the replacement and we may not exclude
large densities of particles close to the origin.  In particular, all
estimates (equivalence of ensembles and local central limit theorems)
need to be uniform over the density. This lack of translation
invariance confines us to one-dimension.

The method presented here may apply to other one-dimensional mean-zero
interacting particle systems. However, instead of replacing a local
function by a function of the empirical measure, one will need to
replace a current multiplied by $N$ by a function of the empirical
measure, as what it is done for non-gradient systems, but without
any space average.

\section{Notation and Results}
\label{sec1}

We consider one-dimensional zero-range processes with periodic
boundary conditions to avoid unnecessary technicalities. This process
is a system of random walks on the discrete torus $\bb T_N = \bb Z / N
\bb Z$ where particles interact infinitesimally only when they are at
the same site. Fix a rate function $g: \bb N_0 = \{0, 1, \dots\} \to
\bb R_+$ with $g(0)=0$, $g(k) >0$, $k\ge 1$, and a finite range
probability measure $p(\cdot)$ on $\bb Z$.  The particle dynamics is
described as follows.  If there are $k$ particles at a site $x$, one
of these particles jumps to site $y$ with an exponential rate
$(g(k)/k) p(y-x)$.

For simplicity, we assume that $p(\cdot)$ is symmetric, but our
results remain true, with straightforward modifications, for any
irreducible, finite-range, mean-zero transition probability $p(\cdot)$.
For the rate function $g$, we assume the next conditions:
\begin{eqnarray*}
& \text{(LG)} & \text{ $ \exists\, a_1 >0$ such that $|g(n+1)-g(n)|
  \leq a_1$ for $n\geq 0$}\;, \\
& \text{(M)} & \text{$\exists\, a_0>0$, $b\geq 1$, such that
$g(n+b)-g(n)>a_0$ for $n\geq 0$}\; .
\end{eqnarray*}
A consequence of (LG), (M) is that $g$ is bounded between two linear
slopes: There is a constant $0<a<\infty$ such that $a^{-1}n \leq g(n)
\leq a n$ for all $n\geq 0$.

Denote by $\Omega_N = \bb N_0^{\bb T_N}$ the state space and by $\xi$
the configurations of $\Omega_N$ so that $\xi(x)$, $x\in \bb T_N$,
stands for the number of particles in site $x$ for the configuration
$\xi$. The zero-range process is a continuous-time Markov chain
$\xi_t$ generated by
\begin{equation}
\label{c0}
(\mc L_N f) (\xi) \;=\; \sum_{x \in \bb T_N} \sum_{z\in\bb Z} p(z) \,
g(\xi(x))\, \big[f(\xi^{x,x+z}) -f(\xi)\big]\; ,
\end{equation}
where $\xi^{x,y}$ represents the configuration obtained from $\xi$ by
displacing a particle from $x$ to $y$:
\begin{equation*}
\xi^{x,y}(z) =
\begin{cases}
\xi(x)-1 & {\rm for \ } z=x \\
\xi(y)+1 &{\rm for \ } z=y \\
\xi(z) &{\rm for \ } z \neq x,y.
\end{cases}
\end{equation*}

Now consider an initial configuration $\xi$ such that $\xi(0) \geq 1$.
Tag one of the particles initially at the origin, and follow its
trajectory $X_t$ jointly with the evolution of the process $\xi_t$.
Specially convenient for our purposes is to consider the process as
seen by the tagged particle defined by $\eta_t(x) = \xi_t(x + X_t)$.
This process is again Markovian, now on the set $\Omega^*_N = \{\eta
\in \Omega_N ; \eta(0) \geq 1 \}$ and generated by the operator $L_N =
L_N^{env} + L_N^{tp}$, where $L_N^{env}$, $L_N^{tp}$ are defined by
\begin{eqnarray*}
(L_N^{env} f) (\eta) &=& \sum_{x \neq 0} \sum_{z\in\bb Z} p(z) \,
g(\eta(x)) \, [f(\eta^{x,x+z})-f(\eta)]\\
&+&  \sum_{z\in\bb Z} p(z) \, g(\eta(0)) \,
\frac{\eta(0) -1}{\eta(0)} \, [f(\eta^{0,z})-f(\eta)]\;,
\end{eqnarray*}
\begin{equation*}
(L_N^{tp} f) (\eta) \;=\; \sum_{z\in\bb Z} p(z) \,
\frac{g(\eta(0))}{\eta(0)} \, [f(\theta_z \eta)-f(\eta)]\;.
\end{equation*}
In this formula, the translation $\theta_z$ is defined by
\begin{equation*}
(\theta_z \eta)(x) =
\begin{cases}
\eta(x+z) & {\rm for \ } x \neq 0,-z \\
\eta(z)+1 &{\rm for \ } x=0 \\
\eta(0)-1 &{\rm for \ } x =-z.\\
\end{cases}
\end{equation*}
The operator $L_N^{tp}$ corresponds to jumps of the tagged particle,
while $L_N^{env}$ corresponds to jumps of the other particles,
called environment.

In order to recover the position of the tagged particle from the
evolution of the process $\eta_t$, let $N_t^z$ be the number of
translations of length $z$ up to time t: $N_{t}^z = N_{t-}^z+1 \iff
\eta_{t} = \theta_z \eta_{t-}$. In this case, $X_t=\sum_z z N_t^z$. As
jumps are not simultaneous, the processes
\begin{equation*}
N_t^z - \int_0^t p(z) \frac{g(\eta_s(0))}{\eta_s(0)} ds
\end{equation*}
are orthogonal martingales and, as $\sum zp(z) = 0$, we see that $X_t$
is a martingale with quadratic variation
\begin{equation*}
\<X\>_t \;=\; \sigma^2 \int_0^t \frac{g(\eta_s(0))}{\eta_s(0)} \, ds\;,
\end{equation*}
where $\sigma^2=\sum_z |z|^2 p(z)$.

We now discuss the invariant measures.  For each $\varphi \geq 0$,
consider the product probability measures $\bar \mu_\varphi = \bar
\mu_\varphi^{N,g}$ in $\Omega_N$ defined by
\begin{equation*}
\bar \mu_\varphi(\xi(x) =k) = \frac{1}{Z(\varphi)}
\frac{\varphi^k}{g(k)!}\; ,
\end{equation*}
where $g(k)! = g(1) \cdots g(k)$ for $k \geq 1$, $g(0)!=1$ and
$Z(\varphi)$ is the normalization constant. $Z(\varphi)$ and $\bar
\mu_\varphi$ are well defined for all $\varphi \geq 0$ due to
conditions $(LG), (M)$. Let $\rho = \rho(\varphi) = \int \eta(0) d
\bar \mu_\varphi$. By conditions (LG), (M), $\varphi \mapsto \rho$ is
a diffeomorphism from $[0,\infty)$ into itself. Define then
$\mu_{\rho}= \bar \mu_{\varphi(\rho)}$, since $\rho$ corresponds to
the density of particles at each site. The measure $\{\mu_\rho : \rho
\ge 0\}$ are invariant for the process $\xi_t$ (cf. \cite{Andjel}).

Due to the inhomogeneity introduced at the origin by the tagged
particle, $\mu_\rho$ is no longer invariant for the process
$\eta_t$. However, a computation shows that the size biased measures
$\nu_\rho$ defined by $d \nu_\rho / d \mu_\rho = \eta(0) / \rho$ are
invariant for the process as seen by the tagged particle, reversible
when $p(\cdot)$ is symmetric.  Here, we take $\nu_0 = \delta_{\mf d_0}$,
the Dirac measure concentrated on the configuration
$\mathfrak d_0$ with exactly one particle at the origin, and note
$\nu_\rho\Rightarrow \delta_{\mf d_0}$ as $\rho\downarrow 0$.

From now on, to avoid uninteresting compactness issues, we define
every process in a finite time interval $[0,T]$, where $T<\infty$ is
fixed. Let $\bb T$ be the unit torus and let $\mc M_+(\bb T)$ be the
set of positive Radon measures in $\bb T$.

Consider the process $\xi_t^N =: \xi_{tN^2}$, generated by $N^2 \mc
L_N$. Define the process $\pi_t^{N,0}$ in $\mc D([0,T],\mc M_+(\bb
T))$ as
\begin{equation*}
\pi_t^{N,0}(du) = \frac{1}{N} \sum_{x \in \bb T_N} \xi_t^N (x)
\delta_{x/N} (du)\;,
\end{equation*}
where $\delta_u$ is the Dirac distribution at point $u$.

For a continuous function $ \rho_0: \bb T \to \bb R_+$, define
$\mu^N_{\rho_0(\cdot)}$ as the product measure in $\Omega_N$ given
by $\mu^N_{\rho_0(\cdot)}(\eta(x)=k)=\mu_{\rho_0(x/N)}(\eta(x)=k)$.
The next result is well known (cf. Chapter V \cite{kl}; see also
\cite{dmp}, \cite{gpv}).

\begin{theorem}
\label{th0}
For each $0\le t\le T$, $\pi_t^{N,0}$ converges in probability to the
deterministic measure $\rho(t,u)du$, where $\rho(t,u)$ is the solution
of the hydrodynamic equation
\begin{equation}
\label{ec0}
\left\{
\begin{array}{l}
\partial_t \rho = \sigma^2 \partial_x^2 \varphi(\rho) \\
\rho(0,u) = \rho_0(u),\\
\end{array}
\right.
\end{equation}
and $\varphi(\rho)= \int g(\xi(0)) d \mu_\rho$.
\end{theorem}

Define now the product measure $\nu^N=\nu_{\rho_0(\cdot)}^N$ in
$\Omega^*_N$
given by $\nu_{\rho_0(\cdot)}^N(\eta(x)=k) = \nu_{\rho_0(x/N)}(\eta(x)
=k)$,
and let $\eta_t^N=:
\eta_{tN^2}$ be the process generated by $N^2 L_N$ and starting from
the initial measure $\nu^N$. Define the empirical measure $\pi_t^N$ in
$\mc D([0,T],\mc M_+(\bb T))$ by
\begin{equation*}
\pi_t^N(du) = \frac{1}{N} \sum_{x\in \bb T_N} \eta_t^N(x) \delta_{x/N}(du).
\end{equation*}
Define also the continuous function $\psi:\bb R_+ \to \bb R_+$ by
$$\psi(\rho)  = \int\big(g(\eta(0))/\eta(0)\big)d\nu_\rho = \left\{\begin{array}{rl} \varphi(\rho)/\rho & \ {\rm
      for \ } \rho>0\\
g(1)& \ {\rm for \ } \rho = 0.\end{array}\right.$$

The next theorems are the main results of this article.  We first
identify the scaling limit of the tagged particle as a diffusion
process:

\begin{theorem}
\label{th2}
Let $x_t^N = X^N_t/N$ be the rescaled position of the tagged particle
for the process $\xi_t^N$. Then, $\{x_t^N : t\in [0,T] \}$ converges
in distribution in the uniform topology to the diffusion $\{x_t : t\in
[0,T]\}$ defined by the stochastic differential equation
\begin{equation}
\label{c9}
d x_t = \sigma \, \sqrt{\psi(\rho(t,x_t))} \, dB_t\; ,
\end{equation}
where $B_t$ is a standard
Brownian motion on $\bb T$.
\end{theorem}

Through this characterization we can describe the evolution of the
empirical measure as seen from the tagged particle:

\begin{theorem}
\label{th1}
$\{\pi_t^N : t\in [0,T]\}$ converges in distribution on $\mc
D([0,T],\mc M_+(\bb T))$ to the measure-valued process $\{\rho(t,u
+x_t)du : t\in [0,T]\}$, where $\rho(t,u)$ is the solution of the
hydrodynamic equation (\ref{ec0}) and $x_t$ is given by \eqref{c9}.
\end{theorem}

Recall that $\eta_0^N$ is distributed according to
$\nu_{\rho_0(\cdot)}^N$.  Denote by $\bb P^N$ the probability measure
in $\mc D([0,T],\Omega_N^*)$ induced by the process $\eta_t^N$, and by
$\bb E^N$ the expectation with respect to this process. Denote also by
$E_\mu[h]$ and $\<h\>_\mu$ the expectation of a function $h: \Omega_N
\to \bb R$ with respect to the measure $\mu$; when $\mu = \nu_\rho$,
let $E_\rho[h]$, $\<h\>_\rho$ stand for $E_{\nu_\rho}[h]$,
$\<h\>_{\nu_\rho}$.  Finally, since in the next sections we consider
only the speeded-up process $\eta_t^N$ we omit hereafter the
superscript $N$.

The plan of the paper is now the following.  After some tightness
estimates in Section \ref{sec3}, certain limits are established in
Theorem \ref{s5} in Section \ref{sec4}--with the aid of ``global''
and ``local'' hydrodynamics results in Sections \ref{sec5} and
\ref{sec6}--which give the main Theorems \ref{th2} and \ref{th1}.

\section{Tightness}
 \label{sec3}

To keep notation simple, in this section we assume the transition
probability $p(\cdot)$ to be nearest neighbor:
\begin{equation*}
p(1) \;=\; p(-1)\;=\; 1/2
\end{equation*}
so that $\sigma^2=1$. Denote by $\mc C(\bb T)$ the space of real
continuous functions on $\bb T$ and by $\mc C^2(\bb T)$ the space of
twice continuously differentiable functions on $\bb T$. For a function
$G$ in $\mc C(\bb T)$, denote by $\pi_t^N(G)$ the integral of $G$ with
respect to $\pi^N_t$:
\begin{equation*}
\pi_t^N(G) \;=\; \int G(u) \pi_t^N(du)\;=\;
\frac 1N \sum_{x\in \bb T_N} G(x/N) \eta^N_t(x)\;.
\end{equation*}

For $T>0$, denote by $\mc D_T = \mc D([0,T], \mc M_+(\bb T) \times \mc
M_+(\bb T) \times \bb T \times \bb R_+)$ the path space of c\`adl\`ag
trajectories endowed with the Skorohod topology. For $N\ge 1$, let
$Q_N$ be the probability measure on $\mc D_T$ induced by the process
$(\pi_t^{N,0}, \pi_t^N, x^N_t, \<x^N\>_t)$, where $\<x^N\>_t$ stands
for the quadratic variation of the martingale $x^N_t$. We prove in
this section that the sequence $\{Q_N : N\ge 1\}$ is tight, which
follows from the tightness of each component of $(\pi_t^{N,0},
\pi_t^N, x^N_t, \<x^N\>_t)$.

Let $Q_N^0$ be the probability measure in $\mc D([0,T],\mc M_+(\bb
T))$ corresponding to the process $\pi_t^{N,0}$. As mentioned in
Theorem \ref{th0}, $Q_N^0$ converges to the Dirac-$\delta$ measure
concentrated on the path $\rho(t,u)du$, where $\rho$ is the solution
of \eqref{ec0}. Hence, the sequence $\{Q_N^0 : N\ge 1\}$ is tight.

On the other hand, as $\mc M_+ (\bb T)$ is a metrizable space under
the dual topology of $\mc C(\bb T)$, to show that $\{\pi^{N} _\cdot
: N\ge 1\}$ is tight, it is enough to prove tightness of the
projections $\{\pi^{N}_\cdot (G) : N\ge 1\}$ for a suitable set of
functions $G$, dense in $\mc C(\bb T)$. For $G$ in $\mc C(\bb T)$,
let $Q_N^G$ be the measure in $\mc D([0,T], \bb R)$ corresponding to
the process $\{\pi_t^N(G) : 0\le t\le T\}$. Tightness of the
sequence $\{Q_N^G : N\ge 1\}$ follows from Aldous's criteria in the
next lemma.

\begin{lemma}
\label{s8}
The sequence $\{Q_N^G: N\ge 1\}$ is tight if
\begin{itemize}
\item[(i)] For every $t \in [0,T]$ and every $\epsilon>0$, there
  exists $M>0$ such that
\begin{equation*}
\sup_N \bb P^N \Big[ \, |\pi_t^N(G)|>M \, \Big] < \epsilon\;.
\end{equation*}

\item[(ii)] Let $\mc T_T$ be the set of stopping times bounded by
  $T$. Then, for every $\epsilon >0$,
\begin{equation*}
\lim_{\gamma \to 0} \limsup_{N \to \infty} \sup_{\tau \in \mc T_T}
\sup_{\theta \leq \gamma}
\bb P^N \Big[ \, |\pi_{\tau+\theta}^N(G) - \pi_\tau^N(G)| > \epsilon
\, \Big ] =0\; .
\end{equation*}
\end{itemize}
\end{lemma}

\begin{lemma}
\label{s7}
The sequence $\{Q^G_N : N \ge 1\}$, $G$ in $\mc C^2(\bb T)$, is tight.
\end{lemma}

\begin{proof}
An elementary computation shows that for each $G$ in $\mc C(\bb T)$,
\begin{equation}
\label{c1}
\begin{split}
M_t^{N,G} &= \pi_t^N(G)- \pi_0^N(G) - \int_0^t \frac{1}{N} \sum_{x\in
        \bb T_N} (\Delta_N G) (x/N) \, g(\eta_s(x)) \, ds \\
&-\int_0^t \frac{g(\eta_s(0))}{\eta_s(0)} \, \pi_s^N(\Delta_N G) \, ds
+ \int_0^t \frac{2}{N} \frac{g(\eta_s(0))}{\eta_s(0)} \,
(\Delta_N G) (0) \, ds
\end{split}
\end{equation}
is a martingale of quadratic variation $\<M^{N,G} \>_t$ given by
\begin{equation*}
\begin{split}
& \<M^{N,G} \>_t = \frac{1}{N^2}\int_0^t
\sum_{\substack{x \in \bb T_N\setminus\{0\}\\ z\in\bb Z}}
p(z)\, g(\eta_s(x)) \, [\nabla_{N,z}G(x/N)]^2 \, ds \\
& \quad +\; \frac{1}{N^2} \int_0^t \sum_{z\in\bb Z} p(z) \, g(\eta_s(0))
\, \frac{\eta_s(0)-1}{\eta_s(0)} \, [(\nabla_{N,z}G)(0)]^2 \, ds\\
&\quad  +\; \int_0^t \sum_{z\in\bb Z} p(z) \, \frac{g(\eta_s(0))}{\eta_s(0)}
\, \Big(\frac{1}{N} \sum_{x\in \bb T_N} (\nabla_{N,z} G) (x/N)
\eta_s(x) - \frac{1}{N} (\nabla_{N,z}G) (0)\Big)^2 \, ds\; .
\end{split}
\end{equation*}
In these formulas, $\nabla_{N,z} G$, $\Delta_N G$ correspond to the
discrete first and second derivatives of $G$:
\begin{eqnarray*}
\!\!\!\!\!\!\!\!\!\!\!\!\!\!\! &&
(\nabla_{N,z} G) (u) \;=\; N[G(u+z/N) -G(u)]\; , \\
\!\!\!\!\!\!\!\!\!\!\!\!\!\!\! && \quad
(\Delta_N G) (u) \;=\; N^2\sum_{z\in\bb Z} p(z)\, [G(u+z/N)-G(u)]\;.
\end{eqnarray*}

Since the rate function $g$ grows at most linearly and since the total
number of particles is preserved by the dynamics,
\begin{equation*}
\<M^{N,G}\>_t = \int_0^t \sum_{z \in \bb Z} p(z)\,
\frac{g(\eta_s(0))}{\eta_s(0)}\, \Big(\frac{1}{N} \sum_{x\in \bb T_N}
\nabla_{N,z}G(x/N) \eta_s(x)\Big)^2 ds + R_t^{N,G}\; ,
\end{equation*}
where $|R_t^{N,G}| \leq C_0 N^{-2} \sum_{x\in\bb T_N} \eta_0(x)$ and
$C_0$ is a finite constant which depends only on $G$, $g$, $p$ and
$T$. In particular, $\bb E^N[|R_t^{N,G}|] \le C_1 N^{-1}$.

Note that in contrast with the martingale associated to empirical measure $\pi^{N,0}_t$, due to
the jumps of the tagged particle, the martingale $M^{N,G}$ \emph{does
  not} vanish in $L^2(\bb P^N)$. In particular, we may not expect the
convergence of the empirical measure $\pi^N_t$ to a deterministic
trajectory.

We are now in a position to prove the lemma. Condition $(i)$ of Lemma
\ref{s8} is a direct consequence of the conservation of the total
number of particles. In order to prove condition $(ii)$, recall the
decomposition \eqref{c1} of $\pi_t^N(G)$ as an integral term plus a
martingale. The martingale term can be estimated by Chebychev's
inequality and the explicit form of its quadratic variation:
\begin{eqnarray*}
\!\!\!\!\!\!\!\!\!\!\!\!\!\! &&
\bb P^N \Big[ \, |M_{\tau+\theta}^{N,G} -M_{\tau}^{N,G}|>\epsilon\, \Big]
\ \leq \ \frac{1}{\epsilon^2} \bb E^N
\Big[ (M_{\tau+\theta}^{N,G})^2 -(M_{\tau}^{N,G})^2 \Big]  \\
\!\!\!\!\!\!\!\!\!\!\!\!\!\! && \qquad\quad
\leq\ \frac{C_1}{\epsilon^2} \, \Vert G'\Vert_\infty^2 \, \bb E^N
\Big[ \int_\tau^{\tau +\theta} \Big(\frac{1}{N}\sum_{x\in \bb T_N}
\eta_s(x)\Big)^2 ds \Big] + \frac{C_1}{\epsilon^2 N}\\
\!\!\!\!\!\!\!\!\!\!\!\!\!\! && \qquad\quad
\leq \ \frac{C_1 \Vert G'\Vert_\infty^2 \, \theta}{\epsilon^2}
E_{\nu^N_{\rho_0(\cdot)}}\Big[ \Big(\frac{1}{N}\sum_{x\in \bb T_N}
\eta(x)\Big)^2\Big] + \frac{C_1}{\epsilon^2 N}
\end{eqnarray*}
which converges to 0 as $N\uparrow \infty$ and $\gamma\downarrow 0$.
The integral term can be estimated in the same way, using again the
conservation of the total number of particles. This proves the
lemma.
\end{proof}

It remains to consider the scaled position of the tagged particle
$x_t^N$ and its quadratic variation. We recall that $x_t^N$ is a
martingale with quadratic variation
\begin{equation}
\label{c10}
\< x^N \>_t = \int_0^t \frac{g(\eta_s^N(0))}{\eta_s^N(0)}
\, ds \;.
\end{equation}

\begin{lemma}
\label{s2}
The process $\{(x^N, \<x^N\>) : N \ge 1\}$ is tight for the uniform
topology.
\end{lemma}

\begin{proof}
We need to show that
\begin{equation}
\label{c4}
\lim_{\epsilon \to 0} \limsup_{N\to\infty} Q_N \big[ \sup_{|t-s|\le
  \epsilon} |x^N_t - x^N_s| >\delta \big] \;=\; 0
\end{equation}
for all $\delta>0$ and a similar statement for the quadratic variation
$\<x^N\>_t$. Recall that $\sup_{k\ge 1} g(k)/k \le a < \infty$ and
consider a symmetric random walk $Z^N_t$ on the discrete torus $\bb
T_N$ with jump rate $a$ and transition probability $p(\cdot)$. We may
couple $Z^N_t$ and $X^N_t$ in such a way that the skeleton chains are
equal, i.e., that the sequence of sites visited by both processes are
the same, and the holding times of $Z^N$ are always less than or equal
to the holding times of $X^N$. In particular,
\begin{equation*}
\sup_{|t-s|\le \epsilon} |x^N_t - x^N_s| \;\le\;
\sup_{|t-s|\le \epsilon} |z^N_t - z^N_s|
\end{equation*}
if $z^N_t = Z^N_t/N$. Therefore, \eqref{c4} follows from the tightness
in the uniform topology of a rescaled symmetric random walk.

Tightness of the quadratic variation $\<x^N\>_t$ in the uniform
topology is an elementary consequence of its explicit expression
\eqref{c10} and the boundedness of $g(k)/k$.
\end{proof}

\section{Limit points and proof of Theorems \ref{th2}, \ref{th1}}
 \label{sec4}

The following, which characterizes certain limit points, is the main
result of this section, which yields Theorems \ref{th2} and \ref{th1}.

\begin{theorem}
\label{s5}
The sequence $Q_N$ converges in the Skorohod topology to the law $Q$
concentrated on trajectories $\{(\pi_t^0, \pi_t, x_t, A_t) : 0\le t\le
T\}$ such that $\pi_t^0(du) = \rho(t,u) du$, where $\rho$ is the
unique weak solution of \eqref{ec0}; $x_t$ is the solution of the
stochastic differential equation \eqref{c9}; $\pi_t(du) = \rho(t,x_t
+u) du$ and $A_t = \sigma^2 \int_0^t \psi(\rho(s,x_s))$ $ds$.
\end{theorem}

\noindent{\bf Proof of Theorems \ref{th2} and \ref{th1}.} As the
limit $x_t$ is concentrated on continuous paths, Theorem \ref{s5}
straightforwardly implies Theorems \ref{th2} and \ref{th1}. \qed
\medskip

The proof of Theorem \ref{s5} is now divided in a sequence of lemmatta.  Denote
by $\{\tau_u : u\in\bb T\}$ the group of translations in $\bb T$
acting on points, functions and measures.

\begin{lemma}
\label{s3}
All limit points $Q$ of the sequence $\{Q_N : N\ge 1\}$ are
concentrated on trajectories $\{(\pi_t^0, \pi_t, x_t, A_t) : 0\le t\le
T\}$ in which $x_t$ is a continuous square integrable martingale.
\end{lemma}

\begin{proof}
Assume, without loss of generality, that $Q_N$ converges to $Q$.
Since, by Lemma \ref{s2}, $\{x^N : N \ge 1\}$ is tight for the uniform
topology, $Q$ is concentrated on continuous paths $x_t$.  In
particular, $x^N_t$ converges in law to $x_t$ for all $0\le t\le T$.

The martingale property is inherited by $x_t$ because $x^N_t$
converges in law to $x_t$ and
\begin{equation*}
\bb E^N \Big[ (x^N_t)^2 \Big] \;=\; \bb E^N \Big[ \sigma^2
\int_0^t \frac{g(\eta_s(0))}{\eta_s(0)} \, ds \Big]\; \le\;
a \sigma^2 t
\end{equation*}
uniformly in $N$. Therefore, $x_t$ is a square integrable
martingale relative to its natural filtration.
\end{proof}

\begin{lemma}
\label{s6}
All limit points $Q$ of the sequence $\{Q_N : N\ge 1\}$ are
concentrated on trajectories $\{(\pi_t^0, \pi_t, x_t, A_t) : 0\le t\le
T\}$ in which $\pi_t^0(du) = \rho(t,u) du$, where $\rho$ is the unique
weak solution of \eqref{ec0}, and $\pi_t(du) = \tau_{x_t} \pi^0_t(du)
= \rho(t,x_t +u) du$.
\end{lemma}

\begin{proof}
Assume, without loss of generality, that $Q_N$ converges to $Q$.  The
first statement follows from Theorem \ref{th0}. On the other hand, by
Lemma \ref{s3} and since $\rho$ is continuous, $Q$ is concentrated on
continuous trajectories $\{(\pi^0_t, x_t) : 0\le t\le T\}$. Hence, all
finite dimensional distributions (f.d.d.) of $(\pi^{0,N}_t, x^N_t)$
(and therefore of $\tau_{x^N_t} \pi_t^{N,0}$) converge to the
f.d.d. of $(\pi^0_t, x_t)$ ($\tau_{x_t} \pi_t^{0} = \rho(t,u+x_t)
du$). Since $\pi_t^N = \tau_{x^N_t} \pi_t^{N,0}$ and since the
f.d.d. characterize a measure on $\mc D([0,T])$, the lemma is proved.
\end{proof}

For $\varepsilon >0$, denote $\iota_\varepsilon
=\varepsilon^{-1}\mb 1\{(0,\varepsilon]\}(u)$ and
$\alpha_\varepsilon =(2\varepsilon)^{-1}
\mb 1\{(-\varepsilon, \varepsilon)\}(u)$.  For $l\geq 1$
and $x\in \Z$, denote by $\eta_s^l (x)$ the mean number of particles
in a cube of length $2l+1$ centered at $x\in {\bb T}_N$ at time $s\geq
0$:
\begin{equation*}
\eta^l_s(x) \ = \ \frac{1}{2l+1} \sum_{|y-x|\leq l} \eta_s(y)\;.
\end{equation*}
When $s=0$, we drop the suffix ``$s$'' for simplicity.

A function $h:\Omega_N\to\bb R$ is said to be local if it depends only
on a finite number of sites. For a local, bounded function
$h:\Omega_N\to\bb R$, denote by $H(\rho)$ and $\bar{h} (\rho)$ its expectations with
respect to $\nu_\rho$ and $\mu_\rho$ respectively. Thus, $H,\bar{h} :\bb R_+ \to \bb R$ are the
functions defined by
\begin{equation}
\label{barH}
  H(\rho) \;=\; E_{\nu_\rho}\big[ h(\eta) \big], \ \ {\rm and \ \ }
\bar{h} (\rho) \;=\; E_{\mu_\rho}\big[ h(\xi) \big]\;.
\end{equation}
Also, define for $l\geq 1$ the local function $H_l:\Omega_N \to
\bb R$ given by
\begin{equation*}
H_l(\eta) \;=\; H(\eta^l(0)).
\end{equation*}
Then, $\bar{H_l}:\bb R_+ \to \bb R$ is the function
$\bar{H_l}(\rho) = E_{\mu_\rho}[H_l]$.

 A local
 function $h:\Omega_N\to\bb R$ is said to be Lipschitz if there exists
 a finite subset $A$ of $\bb Z$ and a finite constant $C_0$ such that
 \begin{equation}
 \label{Lip-def}
 \big\vert h(\xi) - h(\xi') \big\vert \;\le\;
 C_0 \sum_{x\in A} \big\vert \xi(x) - \xi'(x)\big\vert
 \end{equation}
 for all configurations $\xi$, $\xi'$ of $\Omega_N$.

Consider in particular the local function $h_0(\eta(0)) =
g(\eta(0))/\eta(0)$.  It follows from assumptions (LG), (M) that
$h_0(\cdot)$ is a Lipschitz function, bounded above by a finite
constant and below by a strictly positive constant.

We now characterize the quadratic variation of $x_t$.

\begin{lemma}
\label{s9}
All limit points $Q$ of the sequence $\{Q_N : N\ge 1\}$ are
concentrated on trajectories $\{(\pi_t^0, \pi_t, x_t, A_t): 0\le t\le
T\}$ such that
\begin{equation*}
A_t \;=\; \sigma^2 \int_0^t \psi(\rho(s,x_s))\, ds
\end{equation*}
for all $0\le t\le T$.  Moreover, $A_t$ is the quadratic variation of
the martingale $x_t$.
\end{lemma}

\begin{proof}
Assume, without loss of generality, that $Q_N$ converges to $Q$. Since
$\<x^N\>_t$ is tight for the uniform topology by Lemma \ref{s2},
$\<x^N\>_t$ converges to a limit $A_t$ for all $0\le t\le T$. By
Proposition \ref{l2}, with respect to $h(\eta(0)) =
g(\eta(0))/\eta(0)$ and
$H(\rho) = \psi(\rho)$, and since for each $0\le
t\le T$ the map $\pi_\cdot \to \int_0^t ds \int
\iota_\epsilon(x)\bar{\psi_l}(\pi_s(\tau_x \alpha_\varepsilon)) \, dx$
is continuous for the Skorohod topology,
\begin{equation*}
\lim_{l\rightarrow \infty}\lim_{\epsilon\to 0} \lim_{\varepsilon \to
  0}Q\Big[ \, \Big| A_t - \sigma^2 \int_0^t ds \int \iota_\epsilon(x)
  \bar{\psi_l}(\pi_s(\tau_x\alpha_\varepsilon)) \, dx\Big| \;>\;
  \delta\Big] \;=\;0
\end{equation*}
for all $0\le t\le T$ and $\delta>0$. By Lemma \ref{s6}, $\pi_t(du) =
\rho(t,x_t+u) du$. Also, $\rho(s,\cdot)$ is continuous for $0\le s\le
T$, and $\bar{\psi_l}(a) \to \psi(a)$ as $l\uparrow \infty$ by bounded
convergence.  Then, as $\varepsilon\downarrow 0$, $\epsilon \downarrow
0$, and $l\uparrow \infty$, we have a.s.
\begin{equation*}
\int \iota_\epsilon(x) \int_0^t
\bar{\psi_l}(\pi_s(\tau_x\alpha_\varepsilon)) \, dsdx \ \to \
\int_0^t \psi (\rho (s,x_s)) \, ds\;.
\end{equation*}

It remains to show that $A_t$ corresponds to the quadratic variation
of the square integrable martingale $x_t$.  By \cite[Corollary
VI.6.6]{JS}, $\{(x_t^N, \<x^N\>_t) : 0\le t\le T\}$ converges in law
to $\{(x_t, \<x\>_t) : 0\le t\le T\}$. Since by the first part of the
lemma, $\{(x_t^N, \<x^N\>_t) : 0\le t\le T\}$ converges to $\{(x_t,
A_t) : 0\le t\le T\}$, $\<x\>_t = A_t$. This concludes the proof
of the lemma.
\end{proof}

Recall that the quadratic variation $\<x\>_t$ of a martingale $x_t$ is
equal to $x_t^2 - x_0^2 - 2\int_0^t x_s \, dx_s$ and that $\<x\>_t$
can be approximated in $L^2$ by the sequence of Riemannian sums
$\sum_j (x_{t_{j+1}} - x_{t_j})^2$, as the mesh of a partition $\{t_j
: 1\le j\le M\}$ of the interval $[0,t]$ vanishes. In particular, one
can prove directly in our context the identity between $A_t$ and the
quadratic variation $\<x\>_t$.  \medskip

It follows from the characterization of continuous martingales that
$x_t$ is a time-changed Brownian motion:

\begin{corollary}
\label{s4}
The rescaled position of the tagged particle $\{x^N_t : 0\le t\le T\}$
converges in law to the solution of the stochastic differential
equation
\begin{equation*}
dx_t \;=\; \sigma \sqrt{\psi(\rho(t,x_t))}\, dB_t\;,
\end{equation*}
where $B_t$ is a Brownian motion and $\rho$ is the solution of the
differential equation \eqref{ec0}.
\end{corollary}

\noindent{\bf Proof of Theorem \ref{s5}.}  By Section \ref{sec3},
the sequence $Q^N$ is tight. On the other hand, by Lemma \ref{s6}
and Corollary \ref{s4}, the law of the the first and the third
components of the vector $(\pi_t^{0}, \pi_t, x_t, A_t)$ are uniquely
determined. Since, by Lemmatta \ref{s6}, \ref{s9}, the distribution
of the second and fourth components are characterized by the
distribution of $x_t$, and $\rho(t,x)$, the theorem is proved. \qed

\section{Global replacement lemma}
\label{sec5}

In this section, we replace the full empirical average of a local,
bounded and Lipschitz function in terms of its density field. The
proof involves only a few modifications of the standard
hydrodynamics proof of \cite[Lemma V.1.10, Lemma V.5.5]{kl}.

\begin{proposition}[Global replacement]
\label{p1}
Let $r:\Omega_N \to \bb R$ be a local, bounded and Lipschitz function.
Then, for every $\delta>0$,
\begin{eqnarray*}
\limsup_{\varepsilon \to \infty} \limsup_{N \to \infty}
 \bb P^N \Big[
\int_0^T \frac{1}{N} \sum_{x\in {\bb T}_N} \tau_x
\V_{\varepsilon N}(\eta_s) ds \geq \delta \Big] =0,
\end{eqnarray*}
where
\begin{equation*}
\V_l(\eta) = \Big| \frac{1}{2l+1} \sum_{|y| \leq l} \tau_y r(\eta)
-\bar{r}(\eta^l(0)) \Big|, {\rm \ \ and \ \ } \bar{r}(a) =
E_{\mu_a}[r]\;.
\end{equation*}
\end{proposition}

For two measures $\mu$, $\nu$ defined on $\Omega_N$ (or $\Omega_N^*$),
denote by $\mc H(\mu|\nu)$ the entropy of $\mu$ with respect to $\nu$:
\begin{equation*}
\mc H(\mu|\nu) \;=\; \sup_{f} \Big\{ \int f d\mu \;-\; \log
\int e^f d\nu \Big\}\;,
\end{equation*}
where the supremum is carried over all bounded continuous functions
$f$.

A simple computation shows that the initial entropy
$\mc H(\nu^N_{\rho_0(\cdot)}|\nu_\rho)$ is bounded by $C_0 N$ for some
finite constant $C_0$ depending only on $\rho_0(\cdot)$ and $g$.
Let $f_t^N(\eta)$ be the density of $\eta_t$ under $\bb P^N$ with
respect to a reference measure $\nu_\rho$ for $\rho>0$, and let $\hat{f}_t^N(\eta) = t^{-1} \int_0^t
f_s^N(\eta) ds$.  By standard arguments (cf. Section V.2 \cite{kl}),
\begin{equation*}
\mc H_N(\hat{f}_t^N):=\mc H(\hat{f}_t^Nd\nu_\rho|\nu_\rho) \leq C_0 N \quad
{\rm and} \quad \mc D_N(\hat{f}_t^N) := \Big \<\sqrt{\hat{f}_t^N}
(-L_N \sqrt{\hat{f}_t^N}) \Big\>_\rho \leq \frac{C_0}{N}\;.
\end{equation*}
Consequently, by Chebyshev inequality, to prove Proposition
\ref{p1} it is enough to show, for all finite constants $C$, that
\begin{equation*}
\limsup_{\varepsilon \to 0} \limsup_{N \to \infty}
\sup_{\substack{\mc H_N(f) \leq C N \\ {\mathcal D}_N(f) \leq C/N}} \int
\frac{1}{N} \sum_{x\in {\bb T}_N} \tau_x \mc V_{\varepsilon N} (\eta)
f(\eta) d \nu_\rho =0
\end{equation*}
where the supremum is with respect to $\nu_\rho$-densities $f$.  Notice that we may remove from the sum
the integers $x$ close to the origin, say $|x|\le 2 \varepsilon N$,
because $\mc V_{\varepsilon N}$ is bounded. After removing these
sites, we are essentially in the space homogeneous case. Proposition
\ref{p1} follows from the two standard lemmatta below as in the proof
of \cite[Lemma V.1.10]{kl}.

\begin{lemma}[Global 1-block estimate]
\label{g3}
\begin{equation*}
\limsup_{k \to \infty} \limsup_{N \to \infty} \sup_{\substack{\mc H_N(f)
    \leq C N \\ \mc D_N(f) \leq C/N}} \int \frac{1}{N}
\sum_{|x| >  2\varepsilon N} \tau_x \V_k (\eta) f(\eta) d \nu_\rho =0\;.
\end{equation*}
\end{lemma}

\begin{lemma}[Global 2-block estimate]
\label{g4}
\begin{eqnarray*}
&&\limsup_{k\rightarrow \infty}\limsup_{\varepsilon \rightarrow 0}
\limsup_{N\rightarrow \infty} \sup_{\substack{\mc H_N(f) \leq C N \\ \mc
    D_N(f) \leq C/N}} \\
&&\ \ \ \ \ \qquad
\frac{1}{2N\varepsilon +1} \sum_{|y|\leq
  N\varepsilon}
\int \frac{1}{N} \sum_{|x| > 2\varepsilon N} |\eta^k(x+y) -
\eta^k(x)|f(\eta) d\nu_\rho = 0\;.
\end{eqnarray*}
\end{lemma}
We now indicate the proofs of Lemmatta \ref{g3} and \ref{g4} in relation
to \cite[Sections V.4, V.5 ]{kl}.
\vskip .1cm

\noindent
{\it Proofs of Lemmatta \ref{g3} and \ref{g4}.}
To be brief, we discuss only the proof of Lemma \ref{g3} through some
modifications of the argument in \cite[Section V.4]{kl}, as the proof
of Lemma \ref{g4}, using the modifications for Lemma \ref{g3} given
below, is on similar lines to that in \cite[Section V.5]{kl}.

In the first step of the 1-block estimate we cut-off high
densities. We claim that
\begin{equation*}
\limsup_{A\rightarrow \infty}\limsup_{k\rightarrow \infty}
\limsup_{N\rightarrow \infty} \sup_{\mc H_N(f)\leq CN}\int
\frac{1}{N} \sum_{|x| > 2\varepsilon N} \tau_x \mc V_k(\eta)
\mb 1 \{\eta^k(x) >A\} f(\eta) d\nu_\rho = 0\;.
\end{equation*}
Since $\mc V_k$ is bounded, we may replace it by a constant and
estimate the indicator by $A^{-1} \eta^k(x)$. After a summation by
parts, the expression is easily shown to be less than or equal to $C_0
(A N)^{-1} \sum_{x\not = 0} \eta(x)$ for some finite constant $C_0$.
To conclude it remains to follow the proof of \cite[Lemma V.4.1]{kl},
applying the entropy inequality with respect to $\nu_\rho$ and keeping
in mind that the marginals of $\mu_\rho$ and $\nu_\rho$ coincide on
sites $x\not = 0$.

Define now $\V_{k,A}(\eta)=\V_k(\eta) 1\{\eta^k(0)\leq A\}$. By the
previous argument, it is enough to show that for every $A>0$,
\begin{equation}
\label{e4}
\limsup_{k \to \infty} \limsup_{N \to \infty}
\sup_{\mc D_N(f) \leq C/N}
\int \frac{1}{N} \sum_{x\in {\bb T}_N} \tau_x \mc V_{k,A} (\eta)
f(\eta) d\nu_\rho =0\;.
\end{equation}

The proof is analogous to the homogeneous case. Since the origin does
not appear, both the Dirichlet form $\mc D_N$ and the the measure
$\nu_\rho$ coincide with the Dirichlet form of the space homogeneous
zero-range process and the stationary state $\mu_\rho$. In particular,
all estimates needed
involve only the functionals of the space-homogeneous process already
considered in \cite{kl}.  \qed

\section{Local replacement lemma}
\label{sec6}

In this section, we replace a bounded, Lipshitz function supported
at the origin by a function of the empirical density.


\begin{proposition}[Local replacement]
\label{l2}
For any bounded, Lipschitz function $h: \bb N_0 \to \bb R$,
and any $t>0$,
\begin{equation*}
\limsup_{l\rightarrow \infty}\limsup_{\epsilon \to 0} \limsup_{\varepsilon
  \to 0} \limsup_{N \to \infty} \bb E^N \Big[\,
\Big|\int_0^t h(\eta_s(0)) - \frac{1}{\epsilon N}\sum_{x=1}^{\epsilon
  N} \bar{H_l}(\eta^{\varepsilon N}_s(x))\, ds \Big|\, \Big] \;=\; 0 \;,
\end{equation*}
where $H(\rho) = E_{\nu_\rho}[h]$, $H_l(\eta) = H(\eta^l(0))$,
and $\bar{H_l}(\rho) = E_{\mu_\rho}[H_l]$.
\end{proposition}

In the proof of this lemma, there are two
difficulties. The first and the most important one is the
absence of a spatial average, a crucial point in the standard one and two
blocks estimates since it allows a cut-off of large densities and a
reduction to translation-invariant densities in the estimation of the
largest eigenvalue of a local perturbation of the generator of the
process. Without the density cut-off, the equivalence of ensembles,
and therefore the local central limit theorem, has to be proved
uniformly over all densities. Moreover, this absence of space
average confines us to one-dimension.

A second obstacle is the lack of translation invariance of the
stationary state, turning the origin into a special site.  Functions
$h(\eta(0))$ and $h(\eta(x))$, for instance, have different
distributions. In particular, in contrast with the original zero-range
process, the integral $\int \{ g(\eta(0)) - g(\eta(x))\} f d\nu_\rho $
cannot be estimated by the Dirichlet form of $f$.

The proof of Proposition \ref{l2} is divided in several steps. We
start with a spectral gap for the evolution of the environment
restricted to a finite cube. For $l\ge 1$, denote by $\Lambda_l$ a
cube of length $2l+1$ around the origin: $\Lambda_l = \{-l, \dots,
l\}$ and by $L^{env}_{\Lambda_l}$ the restriction of the environment
part of the generator to the cube $\Lambda_l$:
\begin{eqnarray*}
(L_{\Lambda_l}^{env} f) (\eta) &=&
\sum_{\substack{x\in\Lambda_l \\ x \neq 0}}
\sum_{y\in\Lambda_l} p(y-x) \, g(\eta(x)) \, [f(\eta^{x,y})-f(\eta)]\\
&+&  \sum_{z\in\Lambda_l} p(z) \, g(\eta(0)) \,
\frac{\eta(0) -1}{\eta(0)} \, [f(\eta^{0,z})-f(\eta)]\;.
\end{eqnarray*}
We assume above, without loss of generality, that $l$ is larger than
the range of $p(\cdot)$.

Let $\nu^{\Lambda_l}_\rho$ be the measure $\nu_\rho$ restricted to the
set $\Lambda_l$.  For $j \ge 1$, denote by $\Sigma_{\Lambda_l, j}$ the
set of all configurations in $\Lambda_l$ with at least one particle at
the origin and $j$ particles in $\Lambda_l$, and by $\nu_{\Lambda_l,
  j}$ the measure $\nu^{\Lambda_l}_\rho$ conditioned to
$\Sigma_{\Lambda_l, j}$:
\begin{equation}
\label{c11}
\Sigma_{\Lambda_l, j} \;=\; \Big \{\eta\in \bb N_0^{\Lambda_l} :
\eta(0)\ge 1 , \sum_{x\in\Lambda_l} \eta(x) = j \Big \}\;, \quad
\nu_{\Lambda_l, j} (\cdot) \;=\; \nu^{\Lambda_l}_\rho \big( \cdot
\big| \Sigma_{\Lambda_l, j} \big)\;.
\end{equation}
Note that $\nu_{\Lambda_l, j}$ does not depend on the parameter
$\rho$.

\begin{lemma}
\label{s10}
There exists a finite constant $C_0$ such that
\begin{equation*}
\< f ; f\>_{\nu_{\Lambda_l, j}} \;\le\; C_0\, l^2\, \< f \, (-
L_{\Lambda_l}^{env} f) \>_{\nu_{\Lambda_l, j}}
\end{equation*}
for all $j\ge 1$, all $l \ge 1$ and all functions $f$ in
$L^2(\nu_{\Lambda_l, j})$. In this formula, $\< f ;
f\>_{\nu_{\Lambda_l, j}}$ stands for the variance of $f$ with respect
to $\nu_{\Lambda_l, j}$.
\end{lemma}

\begin{proof}
This result follows from the spectral gap of the zero-range process
proved in \cite{LSV}. Since $g(k)/k$ is bounded above and below by
finite strictly positive constants, an elementary computation shows
that
\begin{equation*}
\< f ; f\>_{\nu_{\Lambda_l, j}} \;=\; \inf_c \<(f-c)^2\>_{\nu_{\Lambda_l,j}}
\le\; a^2 \inf_c\< (f'-c)^2\>_{\mu_{\Lambda_l,j-1}}\; = \;
a^2\< f' ; f' \>_{\mu_{\Lambda_l, j-1}}
\end{equation*}
\begin{equation*}
{\rm and} \quad
\< f' (- \mc L_{\Lambda_l} f') \>_{\mu_{\Lambda_l, j-1}}
\;\le\; a^2 \< f (- L_{\Lambda_l}^{env} f) \>_{\nu_{\Lambda_l, j}}
\end{equation*}
provided $0<a^{-1} \le g(k)/k \le a$ for all $k\ge 1$. In this
formula, $\mc L_{\Lambda_l}$ is the generator of the zero-range
process \eqref{c0} restricted to the set $\Lambda_l$, $\mu_{\Lambda_l,
  j-1}$ is the canonical measure associated to the zero-range process
restricted to the set $\Lambda_l$ with $j-1$ particles, and $f'(\eta)
= f(\eta + \mf d_0)$, where $\mf d_0$ is the configuration with
exactly one particle at the origin and summation of configurations is
performed componentwise.
\end{proof}

\subsection{Local one-block estimate}
For $l\ge 1$, define the function $V_l (\eta)$ by
\begin{equation*}
V_l(\eta) = h(\eta(0)) - H (\eta^l(0))\;
\end{equation*}
where we recall $h$ is a bounded, Lipschitz function, and $H(a) =
E_{\nu_a}[h(\eta(0))]$. In this subsection we give the second step
for the proof of Proposition \ref{l2}:
\begin{lemma}[One-block estimate]
\label{l3}
For every $0\le t\le T$,
\begin{equation*}
\limsup_{l \to \infty} \limsup_{N \to \infty} \bb E^N  \Big[\,
\Big|\int_0^t V_l(\eta_s) \, ds \Big|\, \Big ] =0\;.
\end{equation*}
\end{lemma}

\begin{proof}
Since the initial entropy $\mc H(\nu_{\rho_0(\cdot)}^N|\nu_\rho)$ is
bounded by $C_0 N$, by the entropy inequality,
\begin{equation*}
\bb E^N \Big[\, \Big| \int_0^t V_l(\eta_s)
\, ds\Big| \, \Big]
\;\leq\; \frac{C_0}{\gamma} + \frac{1}{\gamma N}
\log \bb E_\rho \Big[\exp\Big\{\gamma N \Big|\int_0^t
V_l(\eta_s) \, ds \Big| \Big\} \Big]\; ,
\end{equation*}
where $\bb E_\rho$ denotes expectation with respect to the process
starting from the invariant measure $\nu_\rho$.  Using the elementary
inequality $e^{|x|} \leq e^x+e^{-x}$, we can get rid of the absolute
value in the previous integral, considering $h$ and $-h$. In this
case, by Feynman-Kac formula, the second term on the right hand side
is bounded by $(\gamma N)^{-1} T \lambda_{N,l}$, where $\lambda_{N,l}$
is the largest eigenvalue of $N^2 L_N + \gamma N V_l$.  Therefore, to
prove the lemma, it is enough to show that $(\gamma N)^{-1}
\lambda_{N,l}$ vanishes, as $N\uparrow\infty$, $l\uparrow\infty$, for
every $\gamma>0$.

By the variational formula for $\lambda_{N,l}$,
\begin{equation}
\label{ec1}
(\gamma N)^{-1} \lambda_{N,l} \;=\; \sup_f \Big\{ \< V_l \, f \>_\rho
- \gamma^{-1} N \< \sqrt{f}(-L_N \sqrt{f}) \>_\rho  \Big\}\;,
\end{equation}
where the supremum is carried over all densities $f$ with respect to
$\nu_\rho$.
Recall that we denote by $L^{env}_{\Lambda_l}$ the restriction of
the environment part of the generator to the cube $\Lambda_l$.  As
the Dirichlet forms satisfy $\<\sqrt{f}(-L^{env}_{\Lambda_l}
\sqrt{f})\>_\rho \leq \<\sqrt{f}(-L_N \sqrt{f})\>_\rho$, we may
bound the previous expression by a similar one where $L_N$ is
replaced by $L^{env}_{\Lambda_l}$.

Denote by $\hat{f}_{l}$ the conditional expectation of $f$ given $\{\eta(z)
: z\in\Lambda_l\}$.  Since $V_l$ depends on the configuration $\eta$
only through $\{\eta(z) : z\in\Lambda_l\}$ and since the Dirichlet
form is convex, the expression inside braces in \eqref{ec1} is less
than or equal to
\begin{equation}
\label{c2}
\int  V_l \, \hat{f}_{l} \, d \nu^{\Lambda_l}_\rho \;
-\; \gamma^{-1} N \int \sqrt{\hat{f}_{l}} \,
(-L^{env}_{\Lambda_l} \sqrt{\hat{f}_{l}} ) \, d \nu^{\Lambda_l}_\rho \;,
\end{equation}
where, as in \eqref{c11}, $\nu^{\Lambda_l}_\rho$ stands for the
restriction of the product measure $\nu_\rho$ to $\bb
N_0^{\Lambda_l}$.

The linear term in this formula is equal to
\begin{equation*}
\sum_{j\ge 1} c_{l,j} (f) \int V_l \, \hat{f}_{l,j}
\, d \nu_{\Lambda_l, j} \;,
\end{equation*}
where $\nu_{\Lambda_l, j}$ is the canonical measure defined in
\eqref{c11} and
\begin{equation*}
c_{l,j} (f) \;=\; \int_{\Sigma_{\Lambda_l, j}} \hat{f}_{l}
\, d \nu^{\Lambda_l}_\rho  \;, \quad
\hat{f}_{l,j} (\eta) \;=\; c_{l,j}(f)^{-1} \,
\nu^{\Lambda_l}_\rho (\Sigma_{\Lambda_l, j})\, \hat{f}_{l} (\eta)\;.
\end{equation*}
The sum starts at $j=1$ because there is always a particle at the
origin.  Note also that $\sum_{j\ge 1} c_{l,j} (f) =1$ and that $
\hat{f}_{l,j} (\cdot)$ is a density with respect to $\nu_{\Lambda_l, j}$.

By the same reasons, the quadratic term of \eqref{c2} can be written
as
\begin{equation*}
\gamma^{-1} N \sum_{j\ge 1} c_{l,j} (f)
\int  \sqrt{\hat{f}_{l,j}} \, (-L^{env}_{\Lambda_l}
\sqrt{\hat{f}_{l,j}}) \, d \nu_{\Lambda_l, j} \;.
\end{equation*}
In view of this decomposition, \eqref{ec1} is bounded above by
\begin{equation*}
\sup_{j\ge 1} \sup_{f}  \Big\{ \int V_l  \, f \, d \nu_{\Lambda_l, j}
\;-\; \gamma^{-1} N  \int  \sqrt{f} \,  (-L^{env}_{\Lambda_l} \sqrt{f})
\, d \nu_{\Lambda_l, j} \Big\}\;,
\end{equation*}
where the second supremum is carried over all densities with respect
to $\nu_{\Lambda_l, j}$.

Recall that $V_l(\eta) = h(\eta(0)) - H (\eta^l(0))$. Let
$\tilde{H}_l(j/2l+1) = \int h(\eta(0)) \, d \nu_{\Lambda_l, j}$. By
Lemma \ref{lclt1} below, we can replace $H (\eta^l(0))$ by
$\tilde{H}_l (\eta^l(0))$ in the previous expression.  Let
$V_{l,j}(\eta) = h(\eta(0)) - \tilde{H}_l (j/2l+1)$ and notice that
$V_{l,j}$ has mean zero with respect to $\nu_{\Lambda_l, j}$ for all
$j\ge 1$. By Lemma \ref{s10}, $L^{env}_{\Lambda_l}$ has a spectral gap
on $\Sigma_{\Lambda_l, j}$ of order $C_0 l^{-2}$, uniformly in $j$. In
particular, since $h$ is bounded, and the inverse spectral gap of
$L^{env}_{\Lambda_l}$ is order $l^2$ uniformly in $j$, by Rayleigh
expansion \cite[Theorem A3.1.1]{kl}, for sufficiently large $N$,
\begin{eqnarray*}
\!\!\!\!\!\!\!\!\!\!\!\!\!\!\! &&
\int V_{l,j}  \, f \, d \nu_{\Lambda_l, j}
\;-\; \gamma^{-1} N  \int  \sqrt{f} \,  (-L^{env}_{\Lambda_l} \sqrt{f})
\, d \nu_{\Lambda_l, j} \\
\!\!\!\!\!\!\!\!\!\!\!\!\!\!\! && \qquad
\le\; \frac{\gamma N^{-1}}{1-2\|V_l\|_{L^\infty}
C_0l^2\gamma N^{-1}} \int V_{l,j} (-L^{env}_{\Lambda_l})^{-1} V_{l,j}
\, d \nu_{\Lambda_l, j}\\
&& \qquad
\le\; 2 \gamma N^{-1} \int V_{l,j} (-L^{env}_{\Lambda_l})^{-1} V_{l,j}
\, d \nu_{\Lambda_l, j}
\end{eqnarray*}
uniformly in $j \ge 1$. By the spectral gap of $L^{env}_{\Lambda_l}$ again,
this expression is less than or equal to
\begin{equation*}
C_0 l^2 \gamma N^{-1} \int V_{l,j}^2 \, d \nu_{\Lambda_l, j}\;\le\;
C_0(h) l^2 \gamma N^{-1}
\end{equation*}
because $h$ is bounded. This proves that \eqref{ec1} vanishes as
$N\uparrow\infty$, $l\uparrow\infty$, and therefore the lemma.
\end{proof}

\begin{lemma}
\label{lclt1}
For bounded, Lipschitz function $h:\bb N_0 \rightarrow \R$,
\begin{equation*}
\limsup_{l\rightarrow \infty} \sup_{k\geq 0} \Big |
E_{\nu_{\Lambda_l,k}}[h(\eta(0))] -
E_{\nu_{k/|\Lambda_l|}}[h(\eta(0))] \Big | \ = \ 0\; .
\end{equation*}
\end{lemma}

\begin{proof}
Fix $\epsilon>0$ and consider $(l,k)$ such that $k/|\Lambda_l| \leq
\epsilon$.  We may subtract $h(1)$ to both expectations. Since $h$ is
Lipschitz, the absolute value appearing in the statement of the lemma
is bounded by
\begin{equation}
\label{c12}
C(h) \Big\{ \int \{\eta(0) - 1\} \, d\nu_{\Lambda_l,k}
\;+\; \int \{\eta (0) -1\} \, d\nu_{k/|\Lambda_l|} \Big\} \;.
\end{equation}
Note that both terms are positive because both measures are
concentrated on configurations with at least one particle at the
origin.  We claim that each term is bounded by $a^2 \epsilon$.

On the one hand, since $\nu_\rho (d\eta) = \{\eta(0)/\rho\}
\mu_\rho(d\eta)$, the second term inside braces is equal to $\rho^{-1}
\int \eta (0) \{ \eta(0) - 1\} \, d\mu_{k/|\Lambda_l|}$, where $\rho
= k/|\Lambda_l|$.  since $k\le a g(k)$, we may replace $\eta(0)$ by $a
g(\eta(0))$ and perform a change of variables $\eta' = \eta - \mf d_0$
to bound the second term in \eqref{c12} by $a \varphi(\rho) \le a^2
\rho$.

On the other hand, by the explicit formula for $\nu_{\Lambda_l,k}$,
the first term in \eqref{c12} is equal to
\begin{equation*}
\sum \eta(0) \{\eta(0)-1\} \prod_{x\in\Lambda_l} \frac 1{g(\eta(x))!}
\;\Big/\; \sum \eta(0)  \prod_{x\in\Lambda_l} \frac
1{g(\eta(x))!} \;,
\end{equation*}
where both sums are performed over $\Sigma_{\Lambda_l,k}$. Replacing
$\eta(0)$ by $a^{\pm 1} g(\eta(0))$ in the numerator and in the
denominator, we obtain that the previous expression is less than or
equal to $a^2 E_{\mu_{\Lambda_l,k-1}} [\eta(0)] \le a^2
k/|\Lambda_l|$. In last formula, $\mu_{\Lambda_l,k}$ is the product
measure $\mu_\rho$ conditioned on the hyperplane
$\Sigma^0_{\Lambda_l,k} = \{\xi \in \bb N_0^{\Lambda_l} :
\sum_{x\in\Lambda_l} \xi(x) = k\}$.

For $(l,k)$ such that $k/|\Lambda_l| \geq \epsilon$, write
\begin{equation*}
E_{\nu_{\Lambda_l,k}}[h(\eta(0))] -
E_{\nu_{k/|\Lambda_l|}}[h(\eta(0))] \ =\ \frac 1{\rho}
\Big\{
E_{\mu_{\Lambda_l,k}}[h'(\eta(0))] -
E_{\mu_{k/|\Lambda_l|}}[h'(\eta(0))] \Big\}\;,
\end{equation*}
where $\rho = k/|\Lambda_l|$ and $h'(j) = h(j) j$.  By Corollary 6.1
(parts a,b) \cite{LSV} the last difference in absolute value is
bounded by $C(h) \epsilon^{-1} l^{-1}$.
\end{proof}

\subsection{Local two-blocks estimate} In this subsection we show how
to go from a box of size $l$ to a box of size $\epsilon N$:

\begin{lemma}[Two-blocks estimate]
\label{l4}
Let $H : \bb R_+ \to\bb R$ be a bounded, Lipschitz function. For every
$t > 0$,
\begin{equation}
\label{ec2}
\limsup_{l \to \infty} \limsup_{\epsilon \to 0} \limsup_{N \to \infty}
\bb E^N \Big[\, \Big| \int_0^t \big\{H(\eta_s^l(0)) -\frac{1}{\epsilon
  N} \sum_{x =1} ^{\epsilon N} H(\eta_s^l(x)) \big\} ds \Big| \, \Big]
= 0.
\end{equation}
\end{lemma}

The proof of this lemma is very similar to the proof of Lemma
\ref{l3}. The expectation in (\ref{ec2}) is bounded by
\begin{equation*}
\frac{1}{\epsilon N} \sum_{x=2l+1} ^{\epsilon N} \bb E^N
\Big[ \, \Big| \int_0^t \big\{H(\eta_s^l(0)) - H(\eta_s^l(x)) \big\} ds
\Big| \, \Big] \;+\; \frac{C(H)(2l+1)}{\epsilon N}\;\cdot
\end{equation*}

Following the proof of the one-block estimate, we see that it is
enough to estimate, uniformly in $2l+1\leq x\leq \epsilon N$, the quantity
\begin{equation*}
\sup_{f} \Big\{ \<V_{l,x} f\>_\rho -N \gamma^{-1} \< \sqrt{f}(-L_N
\sqrt{f})\>_\rho\Big\},
\end{equation*}
where the supremum, as before, is over all density functions $f$ with
$\int f d\nu_\rho=1$ and $V_{l,x}$ is defined by
\begin{equation*}
V_{l,x}(\eta) = H(\eta^l(0)) -H(\eta^l(x)).
\end{equation*}

Notice that the blocks $\Lambda_l$ and $\Lambda_l(x) =: \{
-l+x,\ldots,l+x\}$ are disjoint. Let $L^{env}_{\Lambda_{l,x}}$ be
the restriction of $L^{env}_N$ to the set $\Lambda_{l,x}=\Lambda_l
\cup \Lambda_l(x)$ and define the operator $L_{l,x}$ by
\begin{equation*}
L_{l,x} f(\eta) = L^{env}_{\Lambda_{l,x}} f(\eta) +
g(\eta(l))[f(\eta^{l,x-l}) -f(\eta)] + g(\eta(x-l))[f(\eta^{x-l,l})
-f(\eta)].
\end{equation*}
The operator $L_{l,x}$ corresponds to the environment generator of a
zero-range dynamics on which particles can jump between adjacent sites
on each box, and between endpoints $l$ and $x-l$. Since $x \leq
\epsilon N$, we see, by adding and subtracting at most $\epsilon N$
terms, that
\begin{equation*}
\<f(- L_{l,x} f)\>_\rho \leq (1+ \epsilon N)\<f(-L_N f)\>_\rho .
\end{equation*}

Then, it is enough to prove that
\begin{equation*}
\sup_f \Big\{ \Big \<\big\{H(\eta^l(0))
- H(\zeta^l(0)) \big\} f\Big\>_{\nu_\rho^{\Lambda_l^*}}\  - \ \frac{1}{2 \epsilon\gamma }
\<\sqrt{f}(- L_{l,l} \sqrt{f})\>_{\nu_\rho^{\Lambda_l^*}} \Big\}
\end{equation*}
vanishes as $\epsilon\downarrow 0$ and $l\uparrow \infty$. In this
formula, the state space is $\bb N_0^{\Lambda_l^*}$, where
$\Lambda_l^* = \{-l, \dots, 3l+1\}$, the configurations of this space
are denote by the pair $\beta = (\eta, \zeta)$, where $\eta$ belongs
to $\bb N_0^{\Lambda_l}$ and $\zeta$ belongs to $\bb N_0^{\{l+1,
  \dots, 3l+1\}}$, expectation is taken with respect to the measure
$\nu_\rho^{\Lambda_l^*}$, the projection of $\nu_\rho$ on
$\Lambda_l^*$, $L_{l,l}$ is the generator of the environment restricted
to the set $\Lambda_l^*$:
\begin{eqnarray*}
(L_{l,l} f) (\beta) &=& \sum_{\substack{x \neq 0, x\in
    \Lambda^*_l\\ y\in\Lambda_l^*}} p(y-x) \,
g(\beta(x)) \, [f(\beta^{x,y})-f(\beta)]\\
&+&  \sum_{z\in\bb Z} p(z) \, g(\beta(0)) \,
\frac{\beta(0) -1}{\beta(0)} \, [f(\beta^{0,z})-f(\beta)]\;;
\end{eqnarray*}
and the supremum is carried over all densities $f$ with respect to
$\nu_\rho^{\Lambda_l^*}$.

Following the proof of the one-block Lemma \ref{l3}, we need only
to prove that
\begin{equation*}
\sup_{k\geq 1}\sup_f\Big \{\int \big\{H(\eta^l(0))
- H(\zeta^l(0)) \big\} \, f\, d\nu_{\Lambda_{l}^*,k} -
\frac{1}{ 2 \gamma \epsilon }\int \sqrt{f}(-L_{l,l} \sqrt{f})\,
d\nu_{\Lambda_{l}^*,k}\Big\}
\end{equation*}
vanishes with limits on $\epsilon$ and $l$ where the supremum is on
densities $f$ with respect to the canonical measure
$\nu_{\Lambda_{l}^*,k}(\cdot) = \nu_\rho^{\Lambda_{l}^*}
(\cdot|\Sigma_{\Lambda_{l}^*,k})$, defined similarly as in
(\ref{c11}), where $\Sigma_{\Lambda_{l}^*,k} = \{\beta\in \bb
N_0^{\Lambda_{l}^*}: \beta(0)\ge 1 , \sum_{y\in \Lambda_{l}^*}
\beta(y) = k\}$.

Let $W_l(\beta)=H(\eta^l(0)) - H(\zeta^l(0))$.  By the Rayleigh
expansion \cite[Theorem A3.1.1]{kl}, spectral gap estimate Lemma
\ref{s10} applied to $L_{l,l}$ (which can be thought of as the
environment generator on a block of length $2|\Lambda_l|$) and
boundedness of $H$, for large $N$ and small $\epsilon$,
\begin{eqnarray*}
&&  \int W_{l} \, f\, d\nu_{\Lambda_{l}^*,k} -
\frac{1}{2 \gamma \epsilon}\int \sqrt{f}(- L_{l,l}
\sqrt{f})\, d\nu_{\Lambda_{l}^*,k}\\
&& \qquad \leq \ \int W_{l}\,  d\nu_{\Lambda_{l}^*,k}
\;+\; \frac{2 \gamma \epsilon} {1- C(H) l^2 \gamma \epsilon}
\int W_{l}\big \{(- L_{l,l})^{-1} W_{l}\big\} \, d\nu_{\Lambda_{l}^*,k}\\
&& \qquad \leq \ \int W_{l} \, d\nu_{\Lambda_{l}^*,k} +
C(H) l^2 \gamma \epsilon
\end{eqnarray*}
for some finite constant $C(H)$ depending on $H$.  The last term
vanishes as $\epsilon\downarrow 0$, while the first term vanishes
uniformly in $k$ as $l\uparrow \infty$ by Lemma \ref{lclt2} below.
\qed

\begin{lemma}
\label{lclt2}
For a bounded, Lipschitz function $H:\bb R_+ \rightarrow \R$, we have
that
\begin{equation*}
\limsup_{l \to \infty} \sup_{k \geq 0} \Big| \,
E_{\nu_{\Lambda_{l}^*,k}}\Big [ H(\eta^l(0)) - H(\zeta^l(0))
\Big ] \, \Big | \;=\;0\;.
\end{equation*}
\end{lemma}

\begin{proof}
Fix $\epsilon>0$. Using that $H$ is Lipschitz, we have that
$|H(\eta^l(0)) - H(\zeta^l(0))| \leq C (H) \{\eta^l(0) +
\zeta^l(0)\}$, and so the expectation appearing in the statement of
the lemma is less than or equal to $C (H) E_{\nu_{\Lambda_{l}^*,k}} [
\eta^l(0) + \zeta^l(0) ]$. A computation, similar to the one presented
in the proof of Lemma \ref{lclt1}, shows that
\begin{equation*}
E_{\nu_{\Lambda_{l}^*,k}} [ \beta(0)] \;\le\; a^2 \Big\{ 1+
E_{\mu_{\Lambda_{l}^*,k-1}}[\xi(0)]\Big\}\;,\quad
E_{\nu_{\Lambda_{l}^*,k}} [ \beta(y)] \;\le\; a^2
E_{\mu_{\Lambda_{l}^*,k-1}}[\xi(0)]
\end{equation*}
for all $y\not = 0$. In this formula, $\mu_{\Lambda_{l}^*,k}$ stands
for the canonical measure defined by $\mu_{\Lambda_{l}^*,k} (\cdot) =
\mu_\rho^{\Lambda_{l}^*} (\cdot | \sum_{x\in \Lambda_{l}^*} \beta(x) =
k)$, where $\mu_\rho^{\Lambda_{l}^*}$ is the product measure
$\mu_\rho$ restricted to the set $\Lambda_{l}^*$.  In particular, the
expectation appearing in the statement of the lemma is less than or
equal to $C(H) a^2 \{1 + k\}/2|\Lambda_l|$. This concludes the proof
of the lemma in the case $k/2|\Lambda_l| \leq \epsilon$.

Assume now that $k/2|\Lambda_l|\geq \epsilon$. By definition of the
canonical measure $\nu_{\Lambda_{l}^*,k}$ and the grand-canonical
measure $\nu_\rho^{\Lambda_{l}^*}$, the expectation appearing in the
statement of the lemma is equal to
\begin{equation}
\label{e5}
\frac 1{E_{\mu_{\Lambda_{l}^*,k}} [ \beta(0) ]} \,
E_{\mu_{\Lambda_{l}^*,k}}\Big [ \beta(0) \Big\{ H(\eta^l(0)) -
H(\zeta^l(0)) \Big\} \Big] \;.
\end{equation}
Since the measure is space homogeneous, the denominator is equal to
$\rho_{l,k} = k/2|\Lambda_l|$, while in the numerator we may replace
$\beta(0)$ by $\eta^l(0)$. The numerator can therefore be rewritten as
\begin{equation*}
E_{\mu_{\Lambda_{l}^*,k}}\Big [ \Big\{ \eta^l(0) - \rho_{l,k}\Big\}
\Big\{ H(\eta^l(0)) - H(\zeta^l(0)) \Big\} \Big] +
\rho_{l,k} E_{\mu_{\Lambda_{l}^*,k}}\Big [ H(\eta^l(0)) -
H(\zeta^l(0)) \Big]\;.
\end{equation*}
The second term vanishes because the measure $\mu_{\Lambda_{l}^*,k}$
is space homogeneous, while the first one is absolutely bounded by
$C(H) E_{\mu_{\Lambda_{l}^*,k}} [ \, | \eta^l(0) - \rho_{l,k} |\,
]$. By \cite[Corollary 6.1 (C)]{LSV}, this expression is less than or
equal to
\begin{equation*}
C'(H) E_{\mu_{\rho_{l,k}}^{\Lambda_{l}^*}} \Big[ \, \big | \xi^l(0)
- \rho_{l,k} \big |\, \Big ] \;\le\; C'(H) \, \sigma (\rho_{l,k})
\, l^{-1/2}\;,
\end{equation*}
where $\sigma (\rho)$ stands for the variance of $\xi(0)$ under
$\mu_\rho$. By \cite[(5.2)]{LSV} and since $\varphi(\rho)/\rho$ is
bounded below and above, $\sigma (\rho_{l,k})^2 \le C \rho_{l,k}$.
Therefore, if we recall the denominator in \eqref{e5}, we obtain that
\begin{equation*}
E_{\nu_{\Lambda_{l}^*,k}}\Big [ H(\eta^l(0)) - H(\zeta^l(0))
\Big ]\;\le\; \frac {C(H)}{\sqrt{l\, \rho_{l,k}}}\;,
\end{equation*}
which concludes the proof of the lemma since we assumed the density to
be bounded below by $\epsilon$.
\end{proof}

\subsection{Proof of Proposition \ref{l2}}
Recall $H(\rho) = E_{\nu_\rho}[h]$, $H_l(\eta) = H(\eta^l(0))$, and
$\bar{H_l}(\rho) = E_{\mu_\rho}[H_l]$.  Then, we have that
\begin{eqnarray*}
&&\bb E^N \Big[ \, \Big|\int_0^t \big\{h(\eta_s) -
\frac{1}{\epsilon N}\sum_{x=1}^{\epsilon
  N}\bar{H_l}(\eta_s^{\varepsilon N}(x)) \big\} ds \Big|\, \Big] \\
&& \ \ \ \ \leq \; \bb E^N \Big[ \, \Big|\int_0^t \big\{h(\eta_s)
-H(\eta_s^l(0))\big\}ds \Big|\, \Big] \\
&&\ \ \ \ + \; \bb E^N \Big[ \, \Big|\int_0^t \Big\{H(\eta_s^l(0))
-\frac{1}{\epsilon N} \sum_{x=1}^{\epsilon N}
H(\eta_s^l(x))\Big\}ds \Big|\, \Big] \\
&&\ \ \ \ + \; \bb E^N \Big[ \, \Big|\int_0^t \Big\{\frac{1}{\epsilon N}
\sum_{x=1}^{\epsilon N} \Big(H(\eta_s^l(x))-\bar{H_l}
(\eta_s^{\varepsilon N}(x)\Big) \Big\}ds \Big|\, \Big]\;.
\end{eqnarray*}

As $h$ is bounded, Lipschitz, we have $H$ is bounded, Lipschitz by
Lemma \ref{Lip-lemma} below, and so the first and second terms
vanish by Lemmatta \ref{l3} and \ref{l4}. For the third term, we can
rewrite it as
\begin{equation*}
\bb E^N \Big[ \, \Big|\int_0^t \Big\{\frac{1}{N} \sum_{x\in {\bb T}_N}
\iota_\epsilon(x/N)\Big(H(\eta^l_s(x))-\bar{H_l}(\eta_s^{\varepsilon
  N}(x)\Big) \Big\}ds \Big|\,  \Big]
\end{equation*}
where $\iota_\epsilon(\cdot) = \epsilon^{-1}1\{(0,\epsilon]\}$.  In
fact, as $h$ is bounded, we can replace $\iota_\epsilon$ in the last
expression by a smooth approximation.  Then, for fixed $\epsilon>0$
and $l\geq 1$, treating $H_l(\eta) = H(\eta^l(0))$ as a local
function, which is also bounded, Lipschitz as $H$ is bounded,
Lipschitz, the third term vanishes using Proposition \ref{p1} by
taking $N\uparrow \infty$, and $\varepsilon\downarrow 0$. \qed
\medskip

\begin{lemma}
\label{Lip-lemma} Let $h:\Omega_N\to \bb R$ be a local, Lipschitz
function.  Then, $H:\bb R_+ \rightarrow \bb R$ given by $H(\rho) =
E_{\nu_\rho}[h]$ is also Lipschitz.
\end{lemma}
\noindent{\it Proof.} The proof is similar to
that of Corollary II.3.7 \cite{kl} which shows $\bar{h}(\rho) =
E_{\mu_\rho}[h]$ is Lipschitz. Following the
proof of Corollary II.3.7 \cite{kl},
 it is not difficult to show $\{\nu_\rho: \rho\geq 0\}$ is a
 stochastically increasing family, and for $\rho_1< \rho_2$ that
$$|H(\rho_1) - H(\rho_2)| \ \leq \ C_h \sum_{x\in
 A}|E_{\nu_{\rho_1}}[\eta(x)] - E_{\nu_{\rho_2}}[\eta(x)]|$$
 where $C_h$ is the Lipschitz constant of $h$, and $A\subset \bb Z$ corresponds to the support of $h$.
If $A$ does not contain the origin, the proof is the same as for
 Corollary II.3.7 \cite{kl}.

Otherwise, it is enough to estimate the difference $|E_{\nu_{\rho_1}}[\eta(0)] -
E_{\nu_{\rho_2}}[\eta(0)]|$.  When $0=\rho_2<\rho_1$, the difference equals
 $E_{\mu_{\rho_1}}[\eta(0)(\eta(0)-1)]/\rho_1\leq
 a\varphi(\rho_1)\leq a^2\rho_1$ as $a^{-1}\leq g(k)/k,
 \varphi(\rho)/\rho \leq a$ through (LG), (M), and
 $E_{\mu_\rho}[g(\eta(0))f(\eta)] =
 \varphi(\rho)E_{\mu_\rho}[f(\eta +\mf d_0)]$ where $\mf d_0$ is the
 configuration with exactly one particle at the origin.  When
 $0<\rho_2<\rho_1$, the difference equals
$$|\rho_1^{-1}E_{\mu_{\rho_1}}[\eta(0)^2] -
 \rho^{-1}_2E_{\mu_{\rho_2}}[\eta(0)^2]|
\ \leq \  |\sigma^2(\rho_1)/\rho_1 -
 \sigma^2(\rho_2)/\rho_1|
+ |\rho_1-\rho_2|
$$
where $\sigma^2(\rho) = E_{\mu_\rho}[(\eta(0)-\rho)^2]$.  The
Lipschitz estimate now follows by calculating a uniform bound on the
 derivative
$$\partial_\rho \frac{\sigma^2(\rho)}{\rho} \ =\
\frac{m_3(\rho)}{\rho \sigma^2(\rho)} -
\frac{\sigma^2(\rho)}{\rho^2}$$ where $m_3(\rho) =
E_{\mu_\rho}[(\eta(0)-\rho)^3]$. For $\rho$ large, under assumptions
(LG), (M), this is on order $O({\rho}^{-1/2})$ from Lemma 5.2
\cite{LSV} and bound $a^{-1}\leq \varphi(\rho)/\rho\leq a$; on the
other hand, as $\rho\downarrow 0$, the derivative is also bounded.
\qed

\end{document}